\documentclass[letterpaper, 10 pt, conference]{ieeeconf}

\IEEEoverridecommandlockouts               
\renewcommand{\baselinestretch}{0.98}

\usepackage[T1]{fontenc}
\usepackage{amssymb,amsmath,epsfig}
\usepackage[latin5]{inputenc}
\usepackage{graphicx,times,latexsym}
\usepackage{subfigure}
\usepackage{psfrag,cite}
\usepackage{setspace}
\usepackage{algorithmic, algorithm}
\usepackage{color}
\usepackage{url}

\newtheorem{theorem}{Theorem}
\newtheorem{lemma}{Lemma}

\newtheorem{corollary}{Corollary}
\newtheorem{Proposition}{Proposition}
\newcommand{\oprocendsymbol}{\hbox{$\bullet$}}
\newcommand{\oprocend}{\relax\ifmmode\else\unskip\hfill\fi\oprocendsymbol}

\newcommand{\modulo}{\text{mod}}

\newcommand{\real}{\mathbb{R}}
\newcommand{\integers}{\mathbb{N}}

\newcommand{\floor}[1]{\lfloor #1 \rfloor}
\newcommand{\ceil}[1]{\lceil #1 \rceil}

\allowdisplaybreaks

\graphicspath{{epsfiles/}}

\title{%Implicit information in event-triggered control over digital channels \\  in the presence of system disturbances
Event-triggered stabilization of disturbed linear systems\\ over digital channels}

\author{Mohammad Javad Khojasteh, Mojtaba Hedayatpour, Jorge
  Cort{\'e}s,  Massimo Franceschetti\thanks{M. J. Khojasteh and M. Franceschetti are with the Department of
    Electrical and Computer Engineering of University of California,
    San Diego. M. Hedayatpour is with the Faculty of Engineering \& Applied Science, University of Regina, SK, Canada.  J.~Cort{\'e}s is with the Department of Mechanical and
    Aerospace Engineering, University of California, San Diego.\{mkhojasteh,massimo,cortes\}@ucsd.edu, hedayatm@uregina.ca}}
 
\begin{document}
\maketitle
\begin{abstract}
We present an event-triggered control strategy for stabilizing a scalar, continuous-time,  time-invariant, linear system over a digital communication channel having bounded delay, and in the presence of   bounded system disturbance. We propose an encoding-decoding scheme, and determine   lower bounds on the packet size and on the information transmission rate which are sufficient for stabilization. We show that for small values of the delay,  the timing information \emph{implicit} in the triggering events is enough to stabilize the system with any positive rate. In contrast, when the delay increases beyond a critical threshold, the   timing information alone is not   enough to stabilize the system   and the transmission rate  begins to increase. Finally, large values of the delay require transmission rates higher than what prescribed by the classic \emph{data-rate theorem}. The results  are numerically validated using a linearized model of an inverted pendulum.
\end{abstract}

%\marginJC{Do we really need keywords for a conference paper?}

%\marginJC{Title is a bit long. How about "Event-triggered stabilization of linear systems over digital channels under disturbance"}

\section{Introduction}\label{sec:intro}
Networked control systems (NCS)~\cite{hespanha2007survey}, where the feedback loop is closed over a communication channel,  are a fundamental component of cyber-physical systems (CPS)~\cite{kumar,murray2003future}.  
In this context, data-rate theorems    state that the minimum communication rate  to achieve stabilization is equal to the \emph{entropy rate} of the system, expressed by the sum of the logarithms of the unstable modes.
% Data-rate theorems essentially state that  communication rate available in the feedback loop should be at least as large as the entropy rate of the plant.
Early examples  of data-rate theorems appeared in~\cite{wong1999systems,baillieul1999feedback}. Key later contributions appeared in~\cite{Mitter} and \cite{nair2004stabilizability}. These works consider a ``bit-pipe" communication channel, capable of  noiseless transmission of a finite number of bits per unit time evolution of the system.
%with bounded and unbounded support additive disturbance in the dynamics of the linear plant. 
Extensions to noisy communication channels are considered in~\cite{sahai2006necessity,tatikonda2004control,matveev2009estimation,7541730,sukhavasi2016linear}.
Stabilization over time-varying bit-pipe channels, including the erasure channel as a special case, are studied in~\cite{Paolo,Lorenzo}.
%who also relate information-theoretic
%results of stabilization over communication channels to the corresponding network-theoretic results by considering the effects of the dropout probabilities on stabilization~\cite{sinopoli2004kalman}. 
% The works in~\cite{Paolo,Lorenzo} by considering the stabilization problem over time-varying channel bridged the information-theoretic
% results of stabilization over communication channels to with the corresponding network-theoretic results about the effect of the
% dropout probabilities on stabilization~\cite{sinopoli2004kalman}. 
Additional formulations include stabilization of   systems with random open loop gains over bit-pipe channels~\cite{kostina2016rate},   stabilization of 
switched linear 
%liberzon2014finite
systems~\cite{yang2017feedback}, systems with uncertain parameters~\cite{ranade2015control,kostina2016rate}, multiplicative noise~\cite{ding2016multiplicative,ranade2013non}, optimal control~\cite{tatikonda2004stochastic,kostina2016ratemm,toli,khina2016multi}, and  stabilization using event-triggered strategies~\cite{MJK-PT-JC-MF:16-allerton,tallapragada2016event,li2012stabilizing2,pearson2017control,ling2016bit,linsenmayer2017delay}.

This paper focuses on the  case of stabilization using event-triggered communication strategies.
%communication occurs only when needed, and the main goal is to minimize the number of transmissions, while at the same time achieve the control objective~\cite{Tabuada,WPMHH-KHJ-PT:12}. 
In this context, a key observation made in~\cite{Level} is that if there is no delay in the communication process, there are no system disturbances, and the controller has knowledge of the triggering strategy, then it is possible to stabilize the system with any positive  rate of transmission. This apparently counterintuitive result can be explained by noting that   the act of triggering essentially reveals the state of the system, which can then be perfectly tracked by the controller. Our previous work~\cite{OurJournal1}  quantifies the  information implicit in the timing of the triggering events, as a function of the communication delay and for a given triggering strategy, showing a \emph{phase transition} behavior. When there are no system disturbances and the delay in the communication channel is small enough, a positive rate of transmission is all is needed to achieve exponential stabilization. When the delay in the communication channel is larger than a critical threshold, the implicit information in the act of triggering  is not enough for stabilization, and the transmission rate must increase. These results are compared with a time-triggered implementation subject to delay in~\cite{khojasteh2017time}. 
%presents necessary and sufficient conditions for stabilizing a linear, time-invariant system without additive disturbance. Moreover, the work in~\cite{khojasteh2017time} compares the required information transmission rate for time-triggered and event-triggered control for exponential convergence of state estimation error in the presence of bounded delay in the communication channel.  

The literature, however, has not considered  to what extent the implicit information in the triggering events is still valuable in the presence of system disturbances. These disturbances add an additional degree of uncertainty in the state estimation process, beside the one due to the unknown delay, and their effect should be properly accounted for. With this motivation,  we consider stabilization of a linear, time-invariant system subject to bounded disturbance over a communication channel having a bounded delay. In comparison with~\cite{OurJournal1}, we consider here a weaker notion of stability, requiring the state to be bounded at all times beyond a fixed horizon,  but without imposing exponential convergence guarantees. This allows to simplify the treatment and  to derive a simpler event-triggered control strategy. We design an encoding-decoding scheme for this strategy, and show that when the size of the packet transmitted through the channel at every triggering event is above a certain fixed value, then for small values of the delay  our strategy achieves stabilization using only implicit information and transmitting at a rate arbitrarily close to zero. In contrast, for values of the delay above a given threshold, the transmission rate must increase and eventually surpasses the one prescribed 
%in the same setting 
by the classic data-rate theorem. It follows that for small values of the delay, we can successfully exploit the implicit information in the triggering events and compensate for the presence of system disturbances. On the other hand,   large values of the delay imply that   information has been excessively aged and corrupted by the disturbance,  so that increasingly higher communication rates are required. All  results are numerically validated by implementing our strategy to stabilize an inverted pendulum, linearized about its equilibrium point, over a communication channel. Proofs are omitted for brevity and will appear in full elsewhere.

% The rest of the paper is organized as follows: Section~\ref{sec:setup} presents problem formulation. Control strategy, sufficient condition for stabilization and coding/decoding algorithms are presented in Section~\ref{sec:design}. Implementation of the control strategy and simulation results are presented in Section~\ref{sec:sim}. Finally,  Section~\ref{sec:conc} concludes the paper. 

\subsubsection*{Notation}
Throughout the paper, $\real$ and $\integers$ represent the set of real and natural numbers, respectively. Also, $\log$ and $\ln$ represent base $2$ and natural logarithms, respectively. For a function $f : \real \rightarrow
\real^n$ and $t \in \real$, we let $f(t^+)$ denote the right-hand limit of $f$ at $t$, namely $\lim_{s \rightarrow t^+} f(s)$. In addition, $\floor{x}$ (resp. $\ceil{x}$) denote the nearest integer less (resp. greater) than or equal to $x$. We denote the modulo function by $\modulo(x,y)$, whose value is the remainder after division of $x$ by~$y$. $\text{sign}(x)$ denotes the sign of~$x$.

\section{Problem formulation}\label{sec:setup}
\begin{figure}[h]
	\centering
 \includegraphics[scale=0.5]{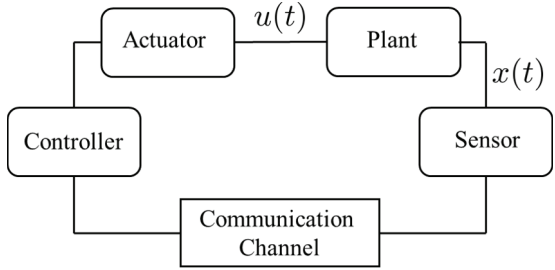}
	\caption{System model.}\label{fig:system} 
\end{figure}

The block diagram of a networked control system as a plant-sensor-channel-controller tuple is represented in Figure~\ref{fig:system}. 
The  plant is described by a scalar, continuous-time, linear time-invariant model as: 
% The plant dynamics are described by a scalar, continuous-time, linear time-invariant (LTI) system 
   \begin{align}\label{syscon}
   	\dot{x}=Ax(t)+Bu(t)+w(t),
   \end{align}
where $x(t) \in \real$ and $u(t) \in \real$ for $t \in [0,\infty)$ are the plant state and control input, respectively, and $w(t)  \in \real$ represents the process disturbance. The latter is upper bounded as: 
\begin{align}\label{noiseupp}
|w(t)|\le M,
\end{align}
where $M$ is a  positive real number. In~\eqref{syscon}, $A$ is a positive real number,  $B \in \real$, and
\begin{align}\label{initialupp}
|x(0)|\le L
\end{align}
for some positive real number~$L$. We assume that the sensor measures the system state exactly, and the controller acts with infinite precision and without delay. However, the measured state is sent to the controller through a communication channel that only supports a finite data rate and is subject to bounded delay. More precisely, when the sensor transmits
packet via the communication channel, the controller will receive the
packet entirely and without any error, but with unknown bounded delay.

The sequence  of triggering times at which the sensor transmits a packet of length $g(t_s^k)$ bits, is denoted by $\{ t_s^k \}_{k \in \integers}$ and the sequence of times at which the controller receives the corresponding packet and decodes it, is denoted by $\{ t_c^k \}_{k \in \integers}$. 
Communication delays are uniformly upper-bounded by $\gamma$, a finite non-negative real number, as follows: 
\begin{align}\label{gammma}
 \Delta_k= t_c^k-t_s^k \leq \gamma ,
 \end{align}
where $\Delta_k$ is the $k^{th}$ \emph{communication delay}. For all $k\geq 1$, we also define the  
 $k^{th}$ \emph{triggering interval} as 
\begin{align}\label{gammma111}
	\Delta'_k = t_s^{k+1}-t_s^{k}.
\end{align}
When referring to a generic triggering or reception time, for  convenience we skip the super-script $k$ in $t_r^k$ and~$t_c^k$.

In this setting, the classical data-rate theorem states that the controller can stabilize the plant if it receives information at least with rate  $A/\ln 2$~\cite{OurJournal1}.
Let $b_s(t)$ be the number of bits transmitted by the sensor up to time $t$. We define the \emph{information transmission rate} as
\begin{align*}
  % \label{Tx-rate}
  R_s = \limsup_{t \rightarrow \infty}\frac{b_s(t)}{t}.
\end{align*}
Since at  every triggering interval the sensor sends $g(t_s)$ bits, we  have  
\begin{align}\label{rss}
  R_s = \limsup_{N \rightarrow \infty} \frac{\sum_{k=1}^{N} g(t_s^k)}{\sum_{k=1}^{N}
    \Delta'_k}.
\end{align}
%The objective is now to precisely quantify the value of $R_s$ required for stabilization when $\gamma$ is in the interval $[0,\infty).$
At the controller, the estimated state is represented by $\hat{x}$ and evolves during the inter-reception times as
%\begin{subequations}\label{eq:ET-strat}
\begin{align}\label{sysest}
  \dot{\hat{x}}(t)=A\hat{x}(t)+Bu(t), \quad t \in [t_c^k,t_c^{k+1}],
\end{align}
starting from $\hat{x}(t_c^{k+})$ with $\hat{x}(0)=\hat{x}_0$. 

We assume that the sensor 
% has causal knowledge of the delay in the communication channel, namely it 
has knowledge of the time the actuator performs the control action. This is to ensure that the sensor can also compute $\hat{x}(t)$ for all time $t$. In practice, this corresponds to assuming an   instantaneous acknowledgment   from the actuator to the sensor via the control input, as discussed in~\cite{sahai2006necessity,ling2017bit}. %This informaiton can be  obtained if the sensor continuously monitors the response time of the plant. Since we need to communicate the reception time to the sensor, 
To obtain such causal knowledge, one can monitor the output of the actuator provided that the control input changes at each reception time.  In case the sensor has only access to the system state, one can use a narrowband signal in the control input to excite a specific frequency of the state,  that can signal the time at which the control action has been applied. 
%As we will discuss later, if the sensor has a causal knowledge of the delay using our event-triggering design it can render $\hat{x}(t)$ from~\eqref{sysest} for all time $t$. 
% We assume that the sensor has a causal knowledge about the delay in the communication channel. In other word, we assume the sensor has the knowledge about the time that actuator gives a new control input to the plant based on the new packet that has been received. This assumption has been showed in Figure~\ref{fig:system} by the dashed line. We will show (in Proposition~\ref{lemmaoffeedback}) that under this assumption with our proposed encoder and decoder scheme the sensor also can calculate $\hat{x}(t)$ at each time $t$. To ensure this assumption we can assume that there is a communication feedback that informs the sensor about the reception time. This communication feedback can also given to the sensor using an acknowledgment in the control input which has been introduce in~\cite{sahai2006necessity} an it has been used in~\cite{ling2017bit} for acknowledge the reception of the packet through
% the control input in for an event-triggering scheme.
The \emph{state estimation error} is defined as
\begin{align}\label{eq:state-estimation-error}
z(t)=x(t)-\hat{x}(t),
\end{align}
where $z(0)=x(0)-\hat{x}_0$. We use this error to determine when a triggering event occurs in our controller design to ensure a property similar to practical stability~\cite{lakshmikantham1990practical} for the system in~\eqref{syscon}.

\section{Control Design}\label{sec:design}
This section proposes our event-triggered control strategy, along with 
 a quantization policy to generate and send  packets at every triggering event, to stabilize the scalar, continuous-time linear time-invariant system described in Section~\ref{sec:setup}. Along the way, we also characterize a sufficient information transmission rate to accomplish this. 
 
Assume a triggering event occurs when
\begin{align}\label{eq:ets}
  |z(t)|=J,
\end{align}
where $J$ is a positive real number. If the controller knows the triggering time $t_s$, then it also knows that $x(t_s)=\pm J+ \hat{x}(t_s)$. It follows that, it may compute the exact value of $x(t_s)$  by just transmitting one single bit at every triggering time.  In general, however, the controller does not have knowledge of $t_s$ because of the delay, but only knows the bound in~\eqref{gammma}. 

Let $\bar{z}(t_c)$ be an estimate of ${z}(t_c)$ constructed by the controller knowing that $|z(t_s)|=v(t_s)$ and using~\eqref{gammma} and the decoded packet received through the communication channel. We define the following updating procedure, called \textit{jump strategy} 
\begin{align}
 \hat{x}(t_c^+)=\bar{z}(t_c)+\hat{x}(t_c).
 \label{eq:jumpst}
\end{align}
At triggering time $t_s$ the sensor encodes the system state in packet $p(t_s)$ of size $g(t_s)$, consisting of the sign of $z(t_s)$ and a quantized version of $t_s$, which we denote by $q(t_s)$, and send it to the controller. Using the bound in \eqref{gammma} and by decoding the received packet, the controller reconstructs the quantized version of $t_s$. Finally, the controller can estimate $z(t_c)$ as follows: 
\begin{align}\label{ctrlapp1}
   \bar{z}(t_c)=\text{sign}(z(t_s)) J e^{A(t_c-q(t_s))}.
\end{align}
Noting that with the jump strategy~\eqref{eq:jumpst}, we have
\begin{align*}
z(t_c^+)=x(t_c)-\hat{x}(t_c^+)=z(t_c)-\bar{z}(t_c),
\end{align*}
the sensor chooses the packet size $g(t_s)$ large enough to 
%performs quantization of the state, and selects a quantization level that
satisfy the following equation for all possible $t_c \in [t_s,t_s+\gamma]$
\begin{align}
|z(t_c^+)| = |z(t_c)-\bar{z}(t_c)| \le \rho_0 J,
  \label{eq:jump-upp}
\end{align}
where $0<\rho_0<1$ is a constant design parameter. To find a lower bound on the size of the packet so that~\eqref{eq:jump-upp} is ensured, the next result bounds how large the difference  $|t_s-q(t_s)|$ of the triggering time and its quantized version can be. 

%\marginJC{I'd always say that having titles for results and remarks helps a lot the reader. We don't have any!}

\begin{lemma}\label{lemmasu}
  For the plant-sensor-channel-controller model with plant dynamics~\eqref{syscon}, estimator dynamics~\eqref{sysest}, triggering strategy~\eqref{eq:ets}, and jump strategy~\eqref{eq:jumpst}, using~\eqref{ctrlapp1} with $J> \frac{M}{A\rho_0}(e^{A\gamma}-1)$, if
   \begin{align}\label{Qua}
   |t_s-q(t_s)| \le \frac{1}{A} \ln(1+\frac{\rho_0-\frac{M}{JA}(e^{A\gamma}-1)}{e^{A\gamma}})
\end{align}
then~\eqref{eq:jump-upp} holds.
\end{lemma}
%The proof of the above result can be found in the %extended version 
%appendix of the paper~\ref{proofoflemma1} available on-line~\cite{online}. 

We next propose our quantization algorithm and rely on  Lemma~\ref{lemmasu} to lower bound the packet size to ensure~\eqref{eq:jump-upp}. 
%The proof of the following result can be found in the 
%appendix of the paper~\ref{proofofthm1}, available on-line~\cite{online}. The design of the quantizer is similar to that our previous work~\cite{OurJournal1}.

\begin{theorem}\label{thm:suf-cond-ET} 
Consider the plant-sensor-channel-controller model with plant dynamics~\eqref{syscon},
  estimator dynamics~\eqref{sysest}, triggering strategy~\eqref{eq:ets}, and
  jump strategy~\eqref{eq:jumpst}. If the control has enough information about $x(0)$ such that state estimation error satisfies  $|z(0)|<J$ , there exists a quantization policy that achieves~\eqref{eq:jump-upp} 
for all $k \in \mathbb{N}$ 
 with a packet size
  \begin{align}\label{Sufi}
  g(t_s^k) \ge \max\left\{0,1+\log \frac{Ab\gamma}{\ln(1+\frac{\rho_0-(M/JA)(e^{A\gamma}-1)}{e^{A\gamma}})}\right\},
  \end{align}
  where $b>1$ and  $J> \frac{M}{A\rho_0}(e^{A\gamma}-1)$.
\end{theorem}

Next, we show that using our encoding and decoding scheme, if the sensor has a causal knowledge of the delay in the communication channel, it can compute the state estimated by the controller.
%The proof of the following result can be found in the appendix of the paper~\ref{proofofpropposition1}
\begin{Proposition}\label{lemmaoffeedback}
  Consider the plant-sensor-channel-controller model with plant dynamics~\eqref{syscon},
  estimator dynamics~\eqref{sysest}, triggering strategy~\eqref{eq:ets}, and
  jump strategy~\eqref{eq:jumpst}.
Using~\eqref{ctrlapp1} and the quantization policy described in Theorem~\ref{thm:suf-cond-ET}, if the sensor has causal knowledge of delay in the communication channel, then the sensor can calculate $\hat{x}(t)$ at each time $t$.
\end{Proposition}

Next, we show that the proposed event-triggered scheme has triggering intervals that are uniformly lower bounded and consequently does not show ``Zeno behavior'', namely infinitely many triggering events in a finite time interval

\begin{lemma}\label{le:mine}
Consider the plant-sensor-channel-controller model with plant dynamics~\eqref{syscon},
  estimator dynamics~\eqref{sysest}, triggering strategy~\eqref{eq:ets}, and
  jump strategy~\eqref{eq:jumpst}. If the packet size satisfies~\eqref{Sufi} for all $k \in \mathbb{N}$, and $J> \frac{M}{A\rho_0}(e^{A\gamma}-1)$ then for all $k \in \mathbb{N}$ 
      \begin{align}\label{lowerinterevent}
t_s^{k+1}-t_s^k \ge \frac{1}{A} \ln (\frac{J+\frac{M}{A}}{\rho_0J+\frac{M}{A}}).
  \end{align} 
\end{lemma}
%The proof can be found at the appendix of the paper~\ref{proofoflemma2}~\cite{online}.
\vspace{3mm}

The frequency with which transmission events are triggered is captured by the triggering rate
\begin{align} \label{trate} R_{tr} &= \limsup_{N\rightarrow
    \infty}\frac{N}{\sum_{k=1}^N \Delta'_k}.
\end{align}
Using Lemma~\ref{le:mine}, we deduce that
\begin{align*}
  R_{tr} \le \frac{A}{\ln (\frac{J+\frac{M}{A}}{\rho_0J+\frac{M}{A}})} 
\end{align*}
for all initial conditions and possible delay and process noise values.
Combining this bound and Theorem~\ref{thm:suf-cond-ET}, we arrive at the following result. 

\begin{figure}[h]
  \centering
  \includegraphics[scale=0.44]{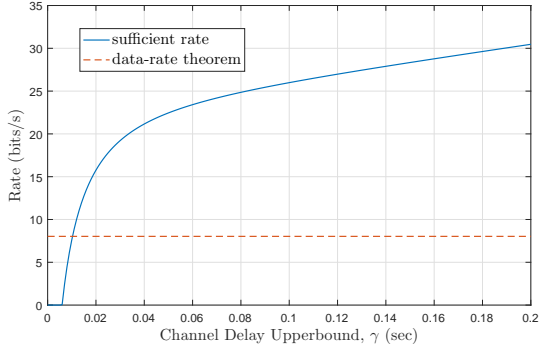} 
  \caption{Illustration of sufficient transmission rate as a function of $\gamma$. Here, A=5.5651, $\rho_0=0.1$, $b=1.0001$, $M=0.2$, and $J=\frac{M}{A\rho_0}(e^{A\gamma}-1)+0.1$.}
  \label{sufficient_rate}
\end{figure}

\begin{corollary}\label{cor:suf-cond-ET2} 
Consider the plant-sensor-channel-controller model with plant dynamics~\eqref{syscon},
  estimator dynamics~\eqref{sysest}, triggering strategy~\eqref{eq:ets}, and
  jump strategy~\eqref{eq:jumpst}. If the control has enough information about $x(0)$ such that state estimation error satisfies  $|z(0)|<J$ with $J> \frac{M}{A\rho_0}(e^{A\gamma}-1)$, there exists a quantization policy that achieves~\eqref{eq:jump-upp} 
for all $k \in \mathbb{N}$ and for all delay and process noise realization
with an information transmission rate
  \begin{align}\label{inftranrate}
& R_s \ge \\\nonumber
 &\frac{A}{\ln (\frac{J+\frac{M}{A}}{\rho_0J+\frac{M}{A}})} \max\left\{0,1+\log \frac{Ab\gamma}{\ln(1+\frac{\rho_0-(M/JA)(e^{A\gamma}-1)}{e^{A\gamma}})}\right\} .
  \end{align}
\end{corollary}
\vspace{3mm}
Figure~\ref{sufficient_rate} shows the sufficient transmission rate as a function of the bound $\gamma$ on the channel delay. As expected, the rate starts from zero and as $\gamma$ increases, goes above the data-rate theorem.

The next result ensures a property similar to practical stability~\cite{lakshmikantham1990practical} for the system in~\eqref{syscon}.
% In the following, we prove that if $|z(0)|\le J$, that is, $|x(0)-\hat{x}(0)| \le J$, then there exit a time $T_0$ and a real number $\kappa$ such that $|x(t)|\le \kappa$ for all $t\ge T_0$.  

\begin{theorem}\label{practicalstb}
Consider the plant-sensor-channel-controller model with plant dynamics~\eqref{syscon},
  estimator dynamics~\eqref{sysest}, triggering strategy~\eqref{eq:ets}, and
  jump strategy~\eqref{eq:jumpst}.
Assume the pair $(A,B)$ is stabilizable.  
  If the control has enough information about $x(0)$ such that state estimation error satisfies  $|z(0)|<J$ with $J> \frac{M}{A\rho_0}(e^{A\gamma}-1)$, and if the sensor use the quantization policy proposed in Theorem~\ref{thm:suf-cond-ET}, then there exists a time $T_0$ and a real number $\kappa$ such that, $|x(t)|\le \kappa$ for all $t\ge T_0$, provided that the packet size is lower bounded by~\eqref{Sufi}.
\end{theorem}
%The proof is available on the on-line appendix of the paper~\ref{practicalstbdd}\cite{online}.

From Corollary~\ref{cor:suf-cond-ET2}, it follows that a transmission rate lower bounded by~\eqref{inftranrate} is sufficient to ensure the property similar to practical stability 
%guarantee 
stated in  Theorem~\ref{practicalstb}.

\section{Simulation}\label{sec:sim}
%\marginJC{This section is way too chatty: most of the discussion can be eliminated with no loss to the paper. Right now, this takes 2.5 pages (not counting plots), which is almost half of the total length we can use!! Specific suggestions below}

%Using the results of Section~\ref{sec:design}, 
%including controller design and coding/decoding algorithms,
We now implement the proposed event-triggered control scheme on a dynamical system such as a linearized inverted pendulum. In this section, initially, a mathematical model of an inverted pendulum mounted on a cart is presented. Then the nonlinear equations are linearized about the equilibrium state of the system. In addition, a canonical transformation is applied to the linear time-invariant system to decouple the equations of motion. %and to determine stable and unstable modes of the system. 

%For each unstable mode of the system (which is one in our example), the proposed control strategy is implemented to stabilize the system. In the end, simulations are used to validate and verify the results. 
%followed by presenting some of the observations and challenges in implementation. 
%\marginJC{Eliminate subsections, it does not make sense for a single example, as we have. Current Section IV.A  can be eliminated entirely, simply write the nonlinear eqs. No need for the plot, the examples is 
%is standard, refer to any book on nonlinear systems for any details regarding modeling.}
%\subsection{Inverted Pendulum}
%We want to validate our result numerically for a linearized model of an inverted pendulum. To do so, 
We consider the two-dimensional problem where motion of the pendulum is constrained in a plane and its position can be measured by angle $\theta$. We assume that inverted pendulum has mass $m_1$, length $l$, and  moment of inertia $I$. Also, the pendulum is mounted on top of a cart of mass $m_2$ constrained to move in $y$ direction. 
%An inverted pendulum of mass $m_1$,Figure~\ref{pendulum}. This is an example of an unstable nonlinear system which has a large number of applications, one of which is balancing a rocket during takeoff. The objective is to balance the pendulum by applying a force $F$ to the cart.  as shown in Figure~\ref{pendulum}. 
%\begin{figure}[t]
%  \centering
%  \includegraphics[scale=0.5]{pendulum.eps} 
%  \caption{A pendulum mounted on a cart.}
%  \label{pendulum}
%\end{figure}
Nonlinear equations governing the motion of the cart and pendulum can be written as follows: 
\begin{align*}
	(m_1+m_2)\ddot{y}+\nu\dot{y}+m_1l\ddot{\theta}\cos\theta-m_1l\dot{\theta}^2sin\theta=F \\
    (I+m_1l^2)\ddot{\theta}+m_1g_0lsin\theta=-m_1l\ddot{y}cos\theta
\end{align*}
where $\nu$ is the damping coefficient between the pendulum and the cart and $g_0$ is the gravitational acceleration. %Additional details  about this example can be found in~\cite{khalil1996noninear}. 

\subsection{Linearizion}
%\marginJC{The English is a bit broken. Also, streamline the exposition, things can be said more succinctly}
We define $\theta=\pi$ as the equilibrium position of the pendulum and $\phi$ as small deviations from $\theta$. We derive the linearized equations of motion using small angle approximation. Let's define state variable $s=[y, \dot{y}, \phi, \dot{\phi}]^T$, where $y$ and $\dot{y}$ are the position and velocity of the cart respectively. Assuming $m_1=0.2$ kg, $m_2=0.5$ kg, $\nu=0.1$ N/m/s, $l=0.3$ m, $I=0.006$ kg/m$^2$, one can write the evolution of $s$ in time as follows:
\begin{align}\label{sysconpendulum}
	\dot{s}=As(t)+Bu(t)+w(t),
\end{align}
\begin{align*}
    \label{AB_matrix}
    A = 
    \begin{bmatrix}
    0 & 1 & 0 & 0 \\
    0 & -0.1818 & 2.6730 & 0 \\
    0 & 0 & 0 & 1 \\
    0 & -0.4545 & 31.1800 & 0
    \end{bmatrix}
    ,  B = 
    \begin{bmatrix}
    0 \\
    1.8180 \\
    0 \\
    4.5450 
    \end{bmatrix}.
\end{align*}

In addition, we add the process noise $w(t)$ to the linearized system model. $w(t)$ is a vector of length four, and  we assume that all the elements of $w(t)$ are upper bounded $M$. 
%as a vector of length four whose elements are upper bounded by $M$, is added to ~\ref{AB_matrix}.
%Also, we added $w(t)$ is the process noise to the linearized system model. $w(t)$ is a vector of length four, and  we assume that all the element of $w(t)$ are upper bounded $M$. 
%and it is assumed to be the same for all four states of the system with appropriate units. 
Also, a simple feedback control law can be derived for ~\eqref{sysconpendulum} as $u=-ks$ where $k$ is chosen such that $A-Bk$ is Hurwitz. We let $k$  be as follows
%\begin{align*}
    $k = \begin{bmatrix}
  -1.00 & -2.04 & 20.36 & 3.93 \end{bmatrix}$.
%\end{align*}

Note that although Theorem~\ref{thm:suf-cond-ET} holds for the linear system with any worst-case delay, the linearizion is only valid for sufficiently small values of $\gamma$. 
% Note that this control law is optimal when the communication channel has infinite capacity and the delay in communication process is zero. There has been efforts for finding optimal feedback control law under communication constraint[Sahai paper, Hassibi and Kostina paper].
%\begin{remark}
%{\rm Since the actual system is nonlinear and the controller developed is linear, it is important as an future problem to investigate the effect of $\gamma$ on controllability of the system. %Hence, our simulation result has been expressed for small values of the $\gamma$. In practice, our result show that for our example when $\gamma < 0.25$ the linearized model is a good approximation of the nonlinear model.
%Results show, for $\gamma\ge0.25$ the linearized model is not a good approximation of the nonlinear model, therefore, for $\gamma\ge0.25$ a nonlinear control strategy will be required to stabilize the system.  
 % }
 % \oprocend
%\end{remark}
%
%\marginJC{The remark can be eliminated, maybe just have 1 sentence in text with the info in the last sentence. }
\subsection{Diagonalization}

The eigenvalues of the open-loop gain of the system $A$ are $e = 
    \begin{bmatrix}
    0 & -5.6041 & -0.1428 & 5.5651
    \end{bmatrix}
$. Hence, three of the four modes of the system are stable and do not need any actuation. Also, the open-loop gain of the system $A$ is diagonalizable (All eigenvalues of $A$ are distinct). As a result, diagonalization of the matrix $A$, enables us to apply Theorem~\ref{thm:suf-cond-ET} to the unstable mode of the system,  and consequently stabilize the whole system.

%\begin{remark} {\rm 
%If the system has more than one unstable modes, we still can use Theorem~\ref{thm:suf-cond-ET} to ensure practical stability. If $\mathfrak{H}$ is the number of unstable modes, we need to append $\mathfrak{h} = \lceil{log_2\mathfrak{H}}\rceil$, number of bits at the beginning of the packet to identify the unstable mode it belongs to.
%In other applications, after diagonalization, if all or some of the resulting SISO systems are unstable, then at every triggering event, we will need to specify which of the SISO systems has been triggered.  
% }
 % \oprocend
%\end{remark}
%Since there is only one positive eigenvalue in matrix $A$, if we could diagonalize the state space representation of the system, then we will end up having four decoupled SISO systems one of which is unstable. This diagonalization, enables us to apply \textit{Theorem 1} to the one unstable SISO system and consequently stabilize the whole system. In the event-triggered control, other stable modes of the system are never triggered and only the one unstable mode will be triggered. 
Using the eigenvector matrix $P$,
% Define \begin{align*}
%     P = 
%     \begin{bmatrix}
%     E(e_1) & E(e_2) & E(e_3) & E(e_4)
%     \end{bmatrix}
% \end{align*}
% where $E(e_i)$ is an eigenvector corresponding to the $i^{th}$ eigenvalue of $A$, then 
% \begin{align*}
%     P = 
%     \begin{bmatrix}
%     1.0000 & 0.0154 & -0.9900 & 0.0147 \\
%     0 & -0.0863 & 0.1414 & 0.0820 \\
%     0 & 0.1750 & 0.0021 & 0.1762 \\
%     0 & -0.9807 & -0.0003 & 0.9808
%     \end{bmatrix}
% \end{align*}
we diagonalize the system to obtain
\begin{equation}\label{AAmj}
    \dot{\tilde{s}} = \tilde{A}\tilde{s}(t)+\tilde{B}\tilde{u}(t)+\tilde{w}(t)
\end{equation}
where
\begin{align*}
    \tilde{A} = 
    \begin{bmatrix}
    0 & 0 & 0 & 0 \\
    0 & -5.6041 & 0 & 0 \\
    0 & 0 & -0.1428 & 0 \\
    0 & 0 & 0 & 5.5651
    \end{bmatrix}
    ,  \tilde{B} = 
    \begin{bmatrix}
    10.0000 \\
    -2.3865 \\
    10.0979 \\
    2.2513 
    \end{bmatrix}
\end{align*}
$\tilde{s}(t) = P^{-1}s(t)$
%$\tilde{A} = P^{-1}AP$, $\tilde{B}=P^{-1}B$,
and $\tilde{w}(t) = P^{-1}w(t)$. Moreover,
$\tilde{u}(t)=-\tilde{k}\tilde{s}(t)$ where $\tilde{k}=kP$, that is, 
%\begin{align*}
    $\tilde{k} = \begin{bmatrix}
    -1.0000 & -0.1295 & 0.7422 & 7.2624
    \end{bmatrix}$.

\subsection{Event-triggered design}
For the first three coordinates of the diagonalized system~\eqref{AAmj} which are stable the state estimation $\hat{s}$ at the controller simply constructs as follows:
\begin{align*}
	%\label{new_sys_est}
    \dot{\hat{s}} = \tilde{A}\hat{s}(t)+\tilde{B}\tilde{u}(t)
\end{align*}
starting from $\hat{s}(0)$.
The unstable mode of the system is as follow
\begin{align}\label{unstablemode}
    \dot{\tilde{s}}_4 = 5.5651\tilde{s}_4(t)+2.2513\tilde{u}(t)+\tilde{w}_4(t)
\end{align}
Then using the problem formulation in section~\ref{sec:setup} the estimated state for the unstable mode $\hat{s}_4$ evolves during the inter-reception times as
\begin{align}\label{esunstablemode}
  \dot{\hat{s}}_4(t)= 5.5651\hat{s}_4(t)+2.2513\tilde{u}(t), \quad t \in [t_c^k,t_c^{k+1}],
\end{align}
starting from $\hat{s}_4(t_c^{k+})$ and $\hat{s}_4(0)$.

The triggering occurs when
\begin{align*}
    |\tilde{z}_4(t)|=|\tilde{s}_4(t)-\hat{s}_4(t)| = J,
\end{align*}
where $|\tilde{z}_4(t)|$ is the estate estimation error for the unstable mode. Let $\lambda_4$ be the eigenvalue corresponding to the unstable mode which is equal to $5.5651$. Then using Theorem~\ref{thm:suf-cond-ET} we choose 
\begin{align*}
 J=\frac{M}{\lambda_4\rho_0}(e^{\lambda_4\gamma}-1)+0.005,
\end{align*}
and 
 the size of the packet for all $t_s$ to be
 \begin{align*}
  g(t_s) = \max\left\{1,\ceil{1+\log \frac{Ab\gamma}{\ln(1+\frac{\rho_0-(M/JA)(e^{A\gamma}-1)}{e^{A\gamma}})}}\right\},
  \end{align*}
  where $b=1.0001$, $\rho_0=0.9$.
  
The packet size for the simulation has two differences from the lower bound provided in Theorem~\ref{thm:suf-cond-ET}. Because the packet size should be an integer we used the ceiling operator, and since we should have at least one bit, to send a packet we take the maximum between 1 and the result of the ceiling operator.

\subsection{Simulation Results}
%\marginJC{This needs a lot of streamlining. We're repeating ourselves at places (e.g., eqs (33) and (34)), and the exposition is a bit wordy.}
The following simulation parameters are chosen for the system: 
simulation time $T=5$ seconds, sampling time $\Delta t=0.005$ seconds, 
%reference set point $r=[0,0,0,0]^T$,
$\tilde{s}(0)={P^{-1}}[0,0,0,0.1001]^T$, and $\hat{s}(0)={P^{-1}}[0,0,0,0.10]^T$.

Theorem~\ref{thm:suf-cond-ET} is developed based on a continuous system but the simulation environments are all digital. We tried to make the discrete model as close to the continuous model by choosing a very small sampling time. However, the minimum upper bound for the channel delay will be equal to one sampling time.
%(in our example, $\Delta t=0.005$ seconds). Hence, in this simulation the worst case delay $\gamma$ is equal to sum of the sampling time and the upper bounded random delay in the communication channel which has been constructed using the random function. 
%This is why in the Fig. \ref{simres} (a), on top, the error is going a little above the triggering function, while in the ideal case, it should be exactly equal to the triggering function.  
%therefore, this problem can be dealt with as if it was like the scalar case. 
A set of three simulations are carried out
%throughout this paper 
as follows. For \textit{simulation (a)} we assumed  the  process disturbance is zero  and channel delay upper bounded by sampling time. In \textit{simulation (b)} we assumed  that the  process disturbance  upper bounded by $M$ and channel delay upper bounded by sampling time. Finally, for \textit{simulation (c)} we assumed  that the  process disturbance  upper bounded by $M$ and channel delay upper bounded by $\gamma$.
%\begin{itemize}
%\item Simulation (a): No process disturbance and channel delay upper bounded by sampling time.
%\item Simulation (b): process disturbance  upper bounded by $M$ and channel delay upper bounded by sampling time.
%\item Simulation (c): process disturbance  upper bounded by $M$ and channel delay upper bounded by $\gamma$
%\end{itemize}
\begin{figure}[t]
  \centering
  \includegraphics[scale=0.44]{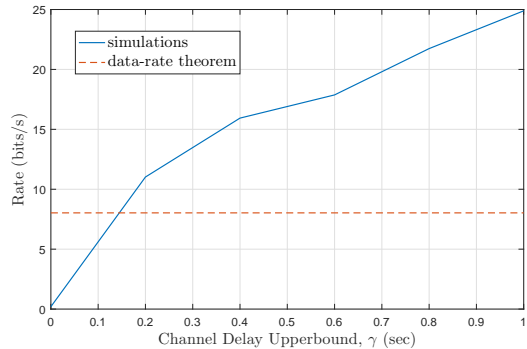} 
  \caption{Information transmission rate in simulations compared to the data-rate theorem. Note that the rate calculated from simulations does not start at zero worst-case delay because the minimum channel delay upper bound is equal to one sampling time (0.005 seconds in this example). $M$ is chosen to be 0.2 in these simulations, and simulation time is $T=5$ seconds. }
  \label{rate}
\end{figure}

\begin{figure*}[t]
\centering
\begin{tabular}{c c c}
	\scriptsize{$M=0.0$, $\gamma=0.005$ sec, $g(t_s)=1$ bit} &
    \scriptsize{$M=0.05$, $\gamma=0.005$ sec, $g(t_s)=1$ bit} &
    \scriptsize{$M=0.05$, $\gamma=0.1$ sec, $g(t_s)=4$ bits} \\
    \includegraphics[width=55mm]{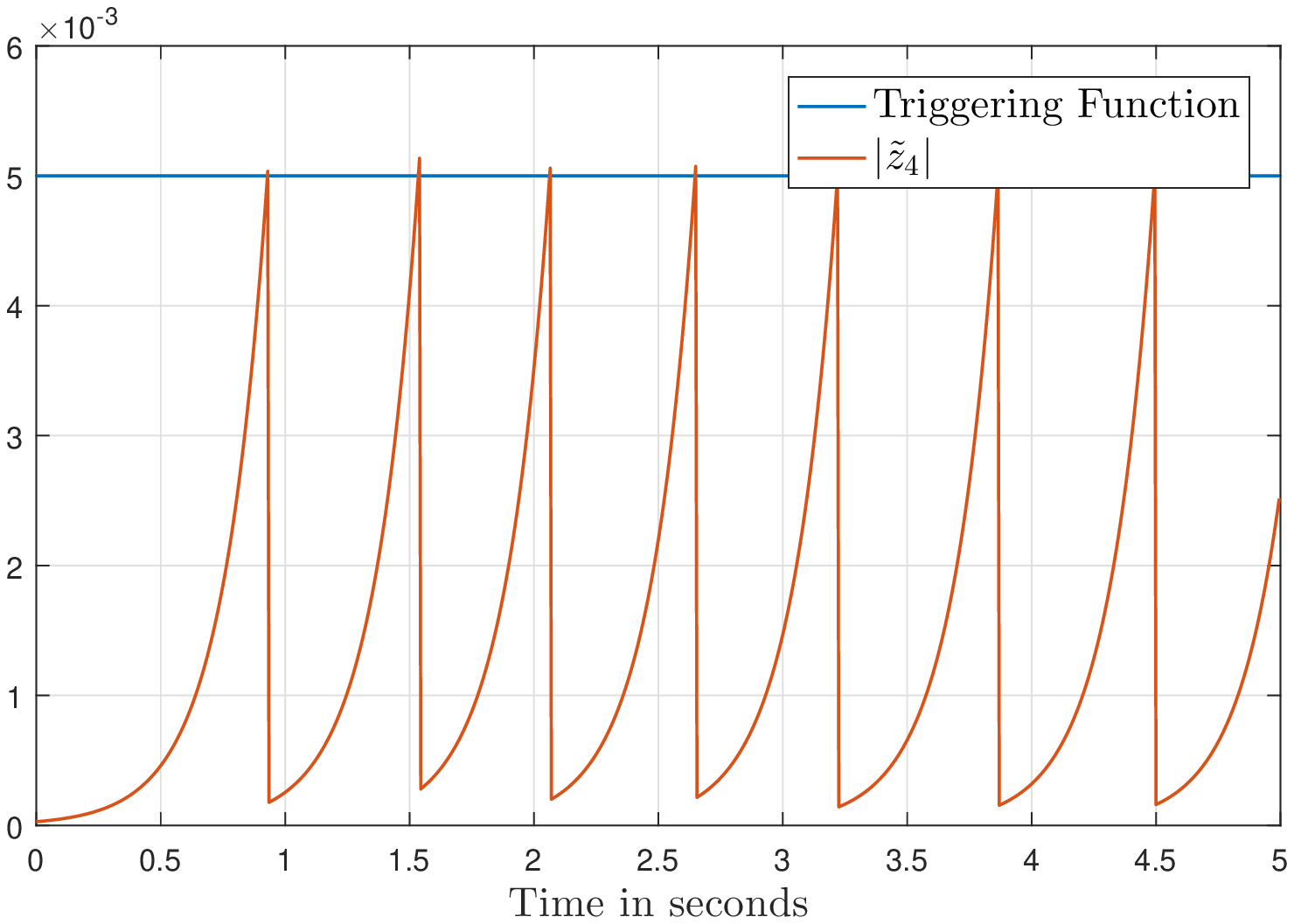} &
    \includegraphics[width=55mm]{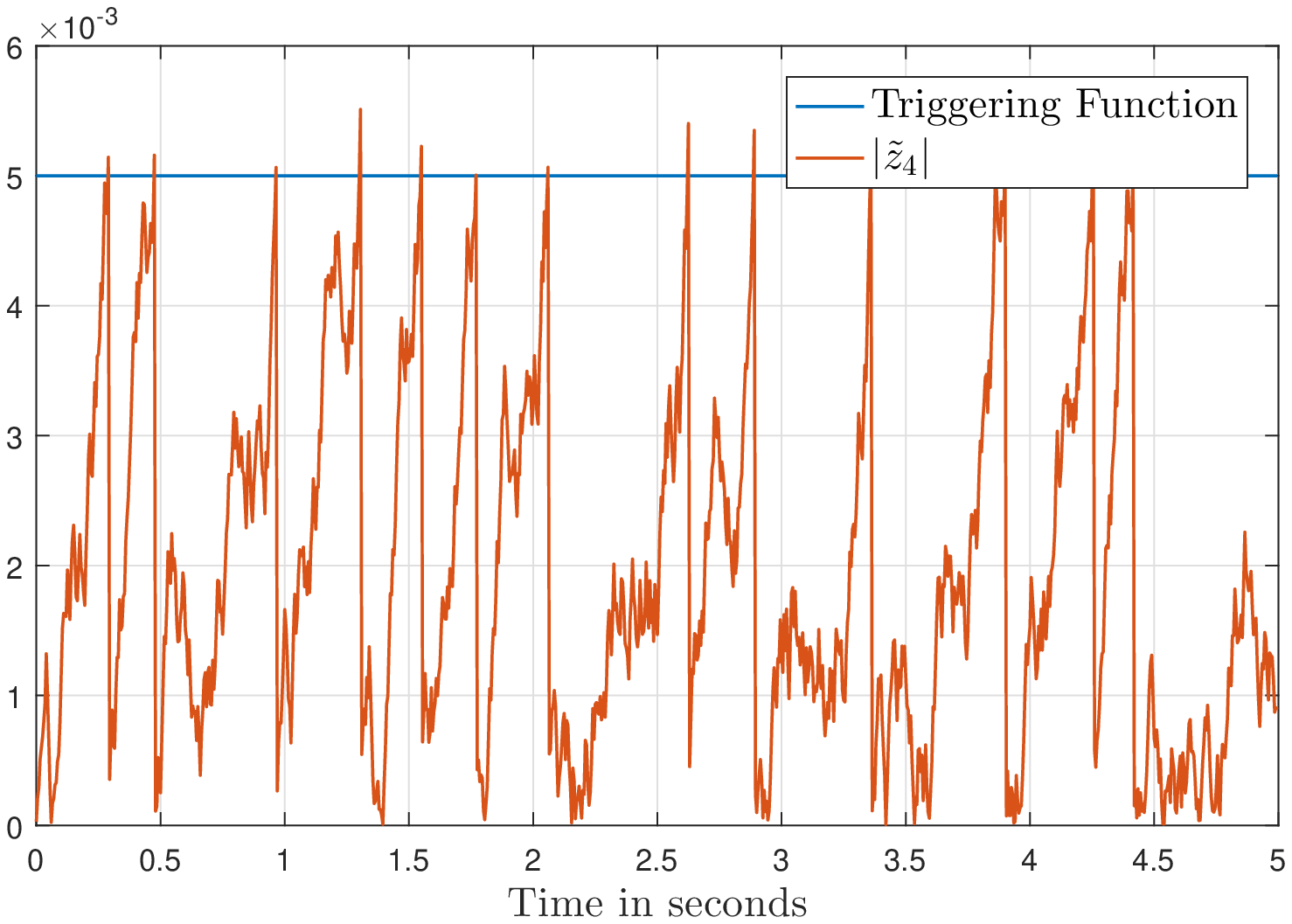} &
    \includegraphics[width=55mm]{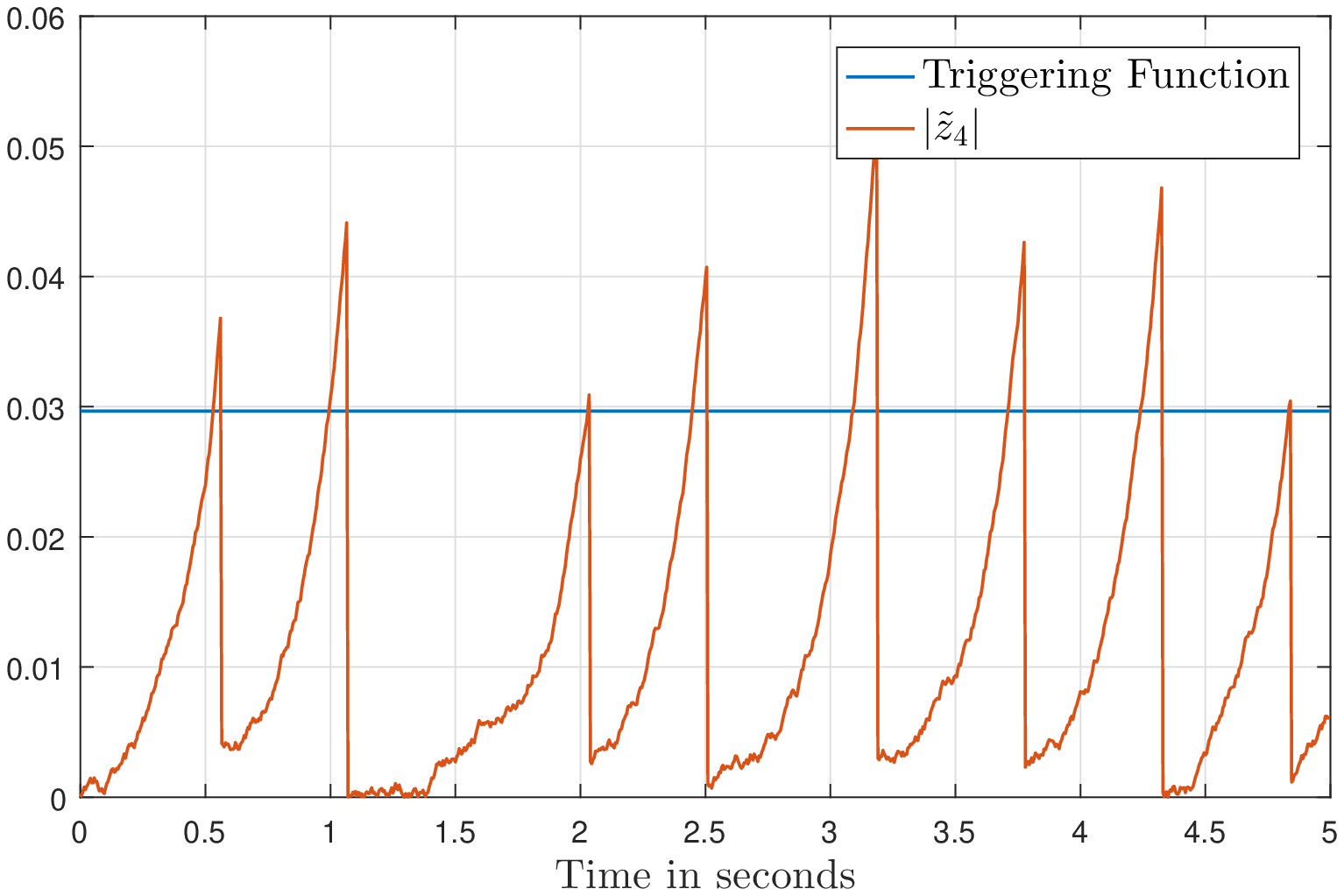} \\
    \includegraphics[width=55mm]{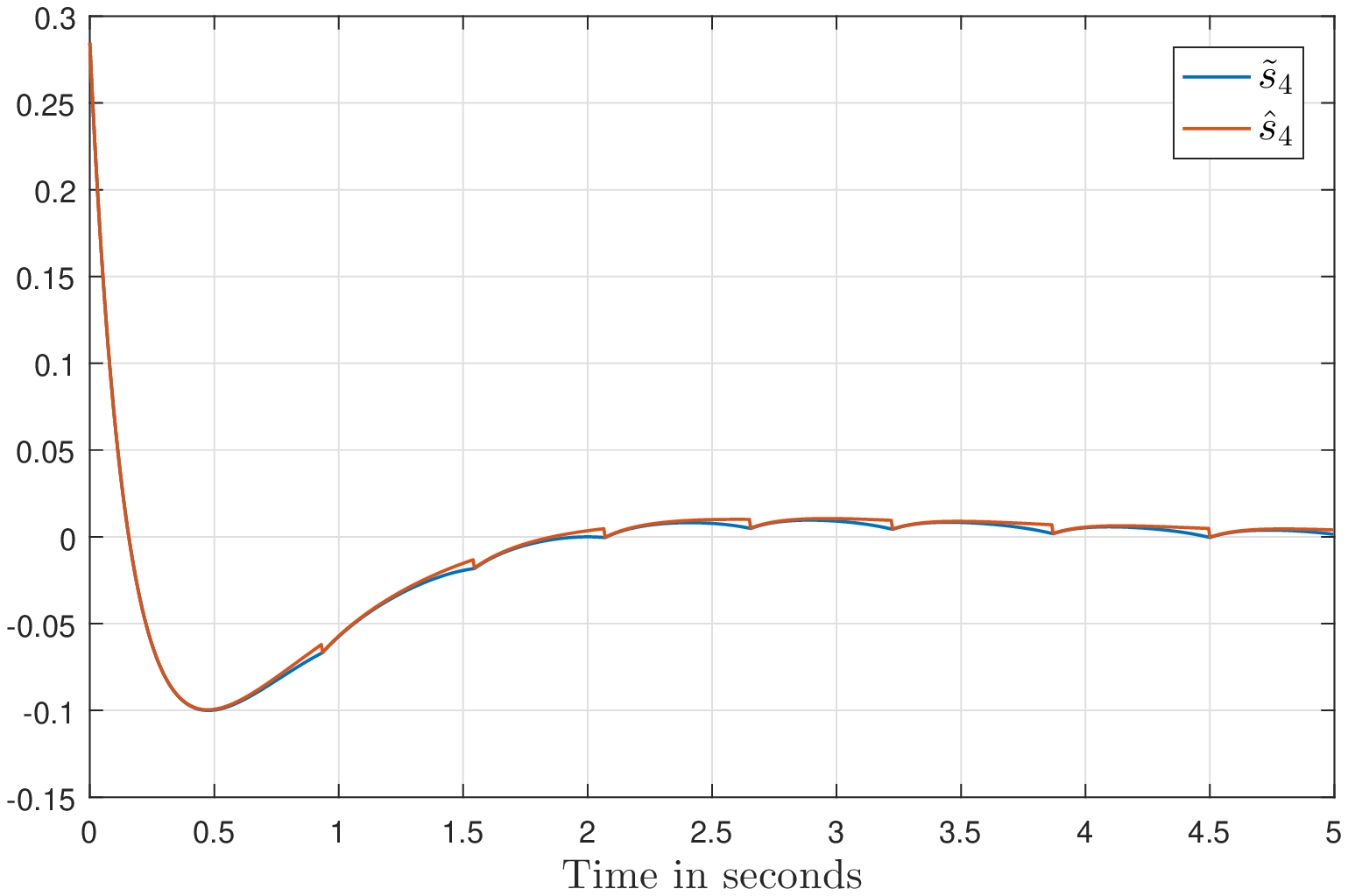} &
    \includegraphics[width=55mm]{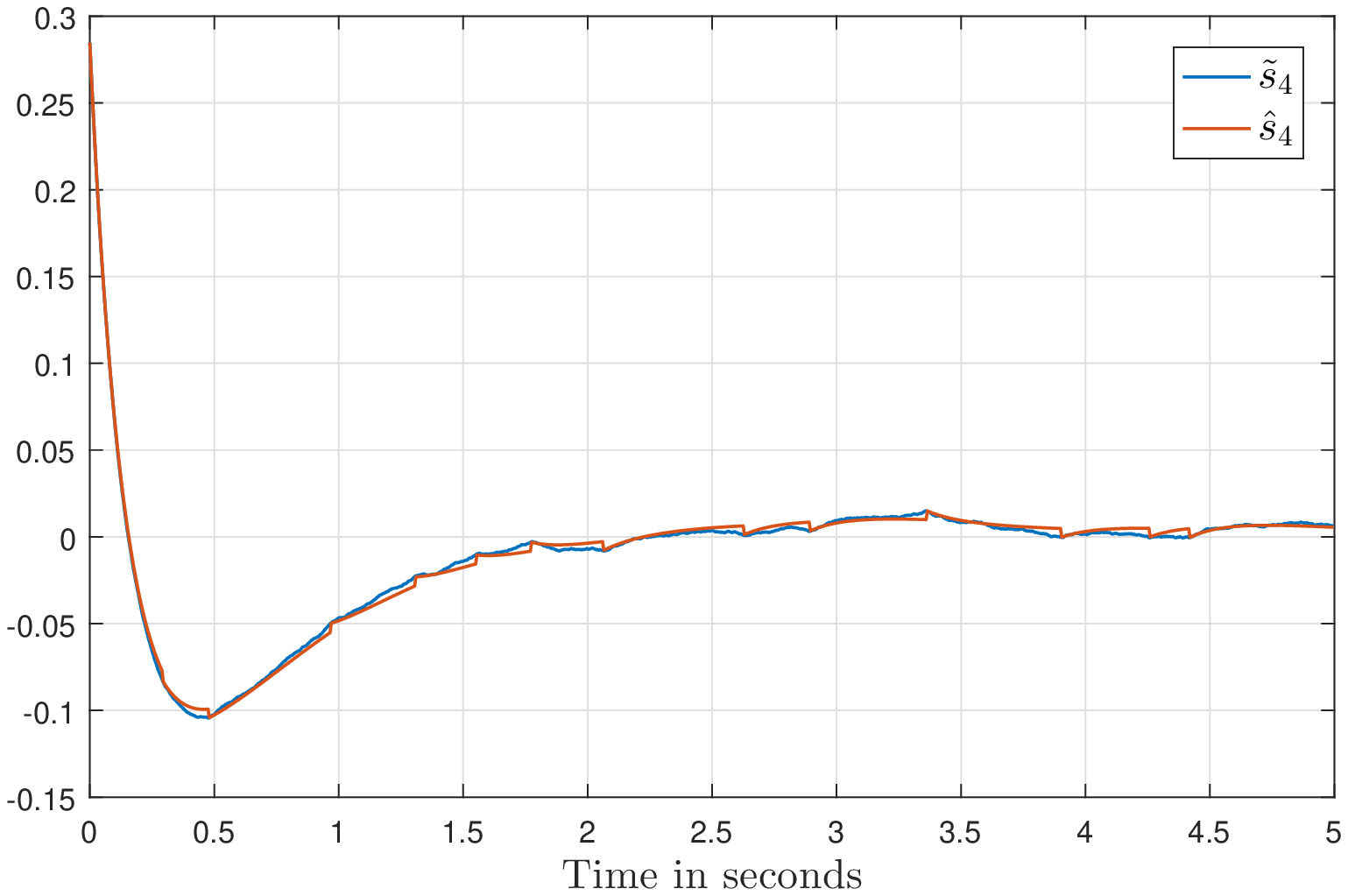} &
    \includegraphics[width=55mm]{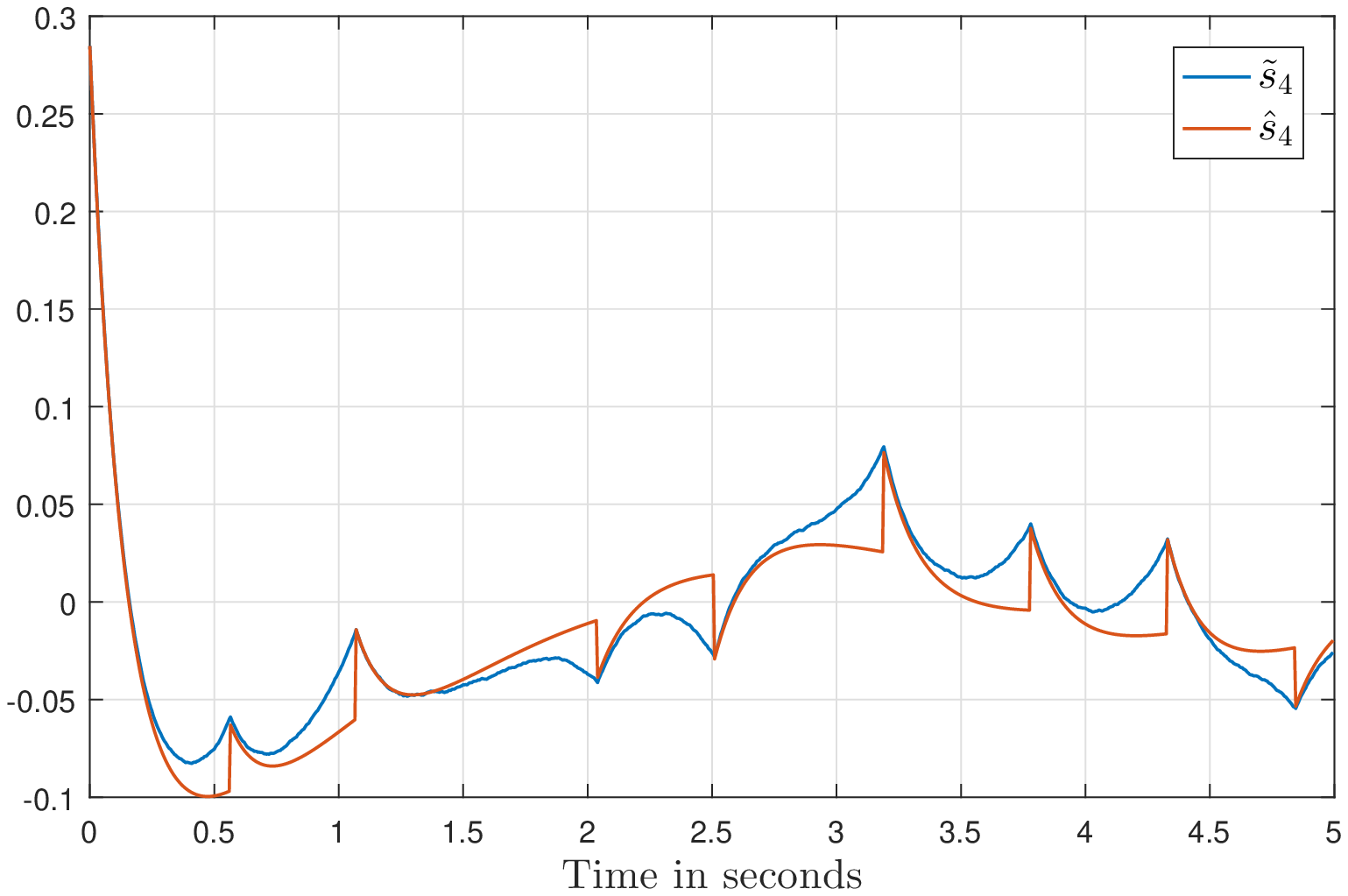} \\
    \includegraphics[width=55mm]{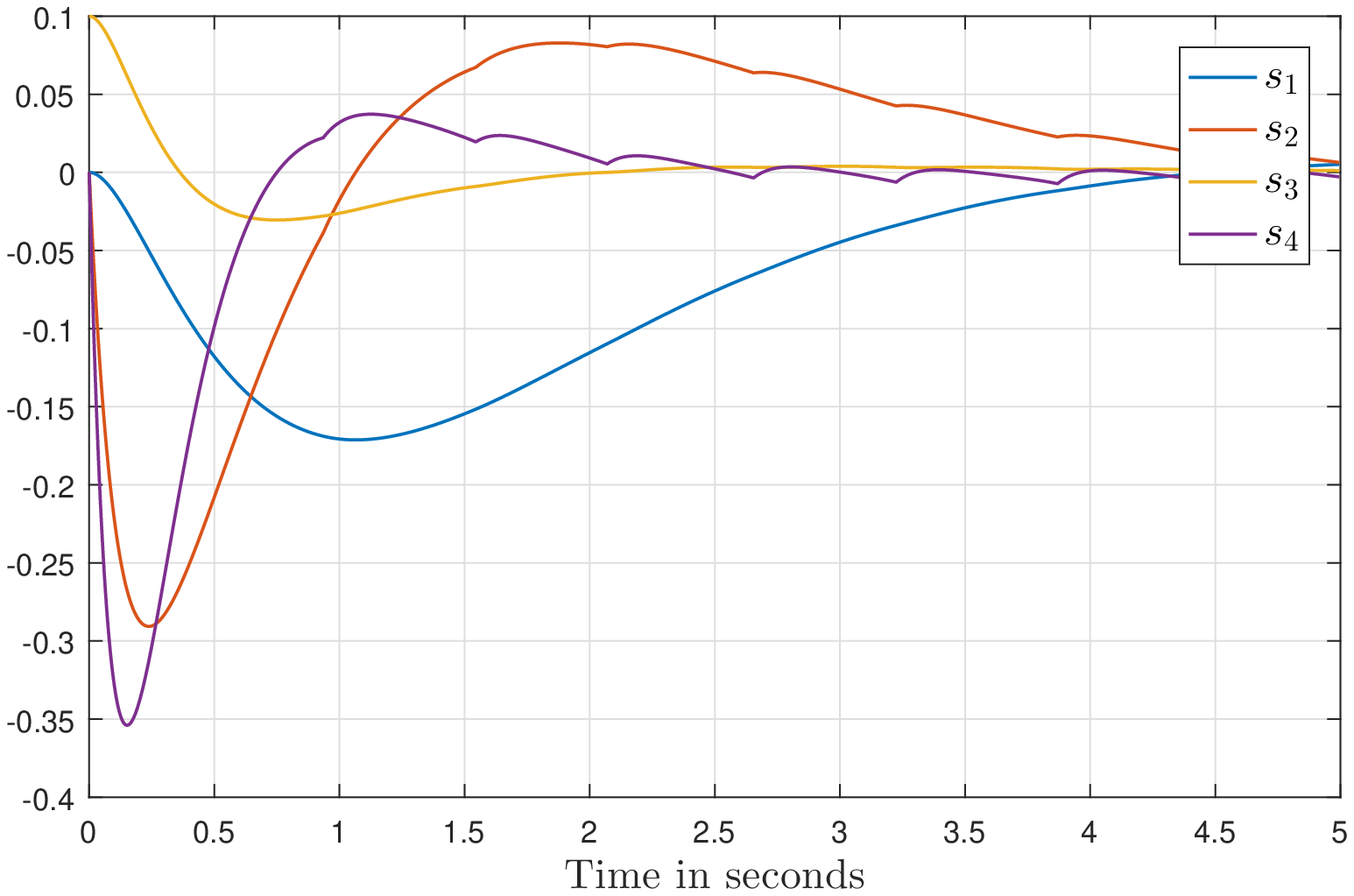} &
    \includegraphics[width=55mm]{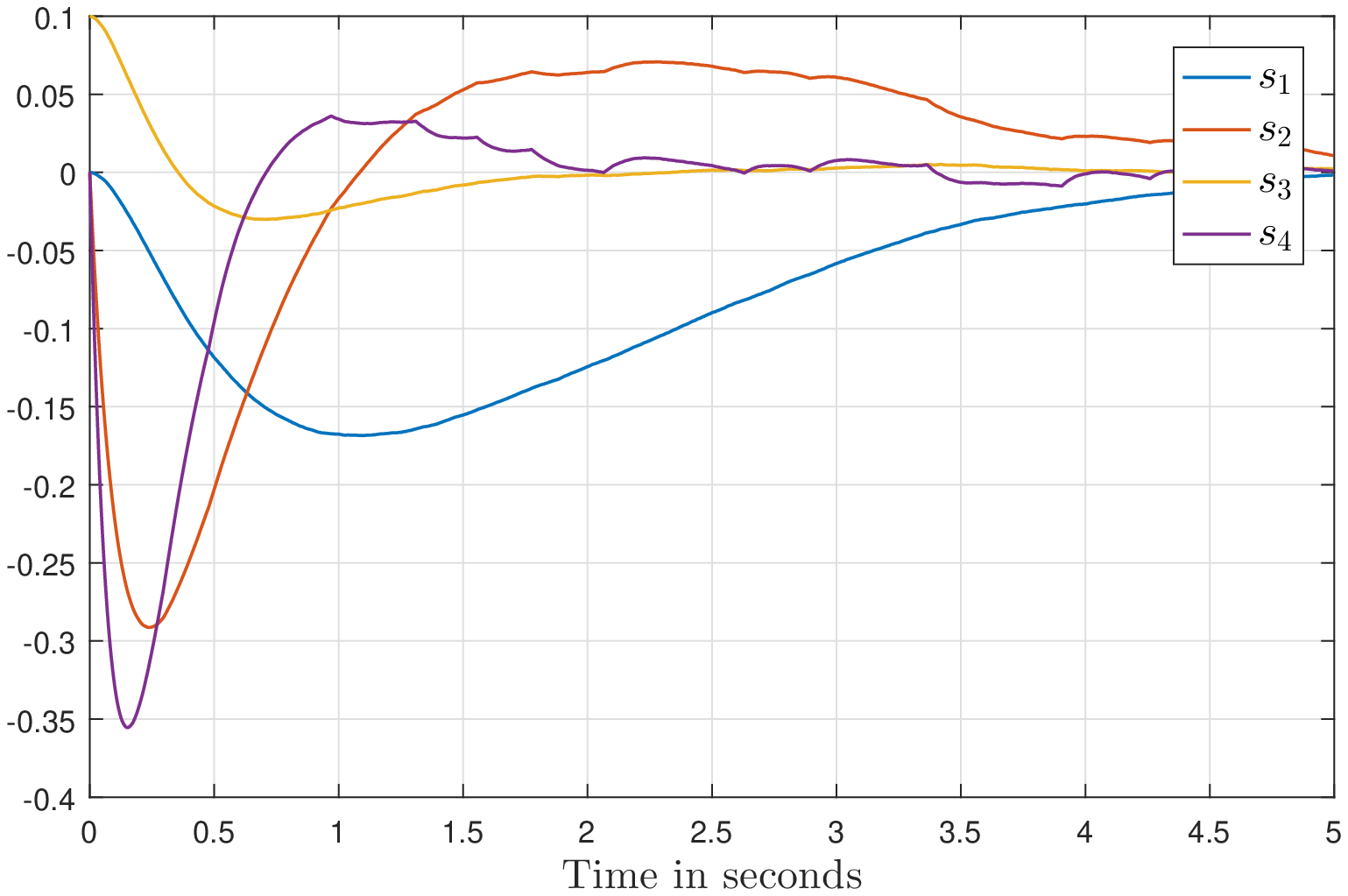} &
    \includegraphics[width=55mm]{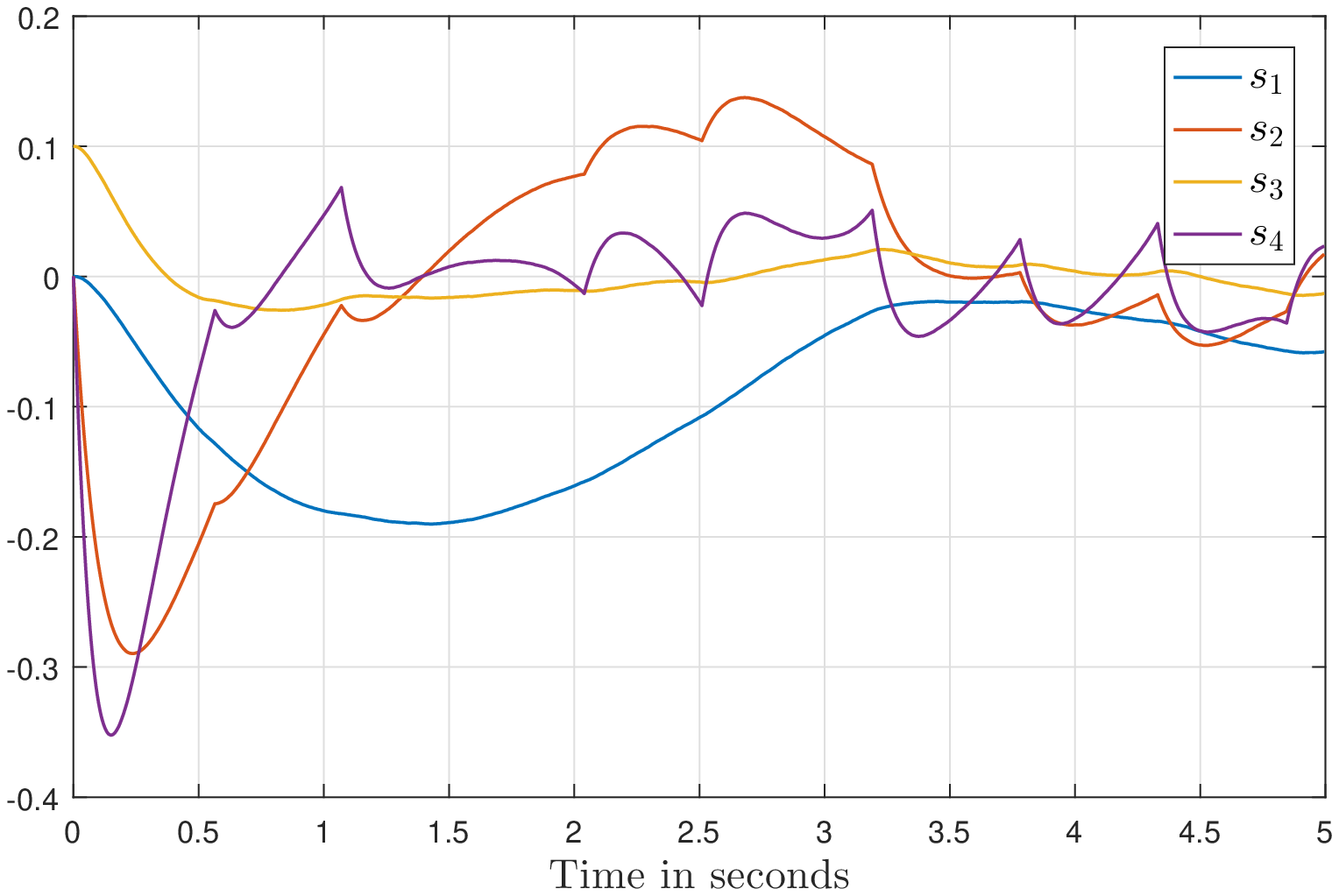} \\
    \small{(a)} & \small{(b)} & \small{(c)}
\end{tabular}
\caption{Simulation results: The first row represents the absolute value of state estimation error for the unstable mode of the system~\eqref{unstablemode}. The second row represents the unstable mode~\eqref{unstablemode}, and its state estimate~\eqref{esunstablemode}. Finally, and the last row represents the evolution of all actual states of the real system~\eqref{sysconpendulum} in time.
%, where the yellow graph is the pendulum angle.
}\label{simres}
\end{figure*}

Simulation results for simulation (a), (b) and (c) are presented in Figure~\ref{simres}. Each column represents a different simulation. The first row shows the triggering function for $\tilde{s}_4$~\eqref{unstablemode} and the absolute value of state estimation error for the unstable coordinate, that is, $|\tilde{z}_4(t)|=|\tilde{s}_4(t)-\hat{s}_4(t)|$. As soon as the absolute value of this error is equal or greater than the triggering function, sensor transmit a packet, and the jumping strategy adjusts $\hat{s}_4$ at the reception time to practically stabilize the system. The amount this error exceeds the triggering function depends on the random channel delay with upper bound $\gamma$. 
In the second row of Figure~\ref{simres}, the evolution of the unstable state~\eqref{unstablemode} and its state estimation are presented~\eqref{esunstablemode}. Finally, the last row in Figure~\ref{simres} represents the evolution of all actual states of the  linearized system~\eqref{sysconpendulum} in time.

Finally, Figure~\ref{rate} presents the simulation of information transmission rate versus the worst-case delay in communication channel $\gamma$ for stabilizing the linearized model of the inverted pendulum.

\section{Conclusions}\label{sec:conc}
We have presented an event-triggered control scheme for the stabilization of noisy, scalar, continuous, linear time-invariant systems over a communication channel subject to random bounded delay. 
We have also developed an algorithm for coding/decoding the quantized version of the estimated states, leading to the characterization of a  
sufficient transmission rate for stabilizing the system. We have illustrated our results on a linearization of the inverted pendulum for different channel delay bounds. Future work  will study the identification of necessary conditions on the transmission rate, the investigation of the effect of delay on nonlinear systems, and the implementation of the proposed control strategies on real systems.

%\marginJC{Check the list of references, some of them are incomplete. For instace, almost no info in [31]}
\section*{Acknowledgements}
This research was partially supported by NSF award CNS-1446891.
\bibliography{mybib} 
\bibliographystyle{IEEEtran}

\end{document}